\documentclass[a4paper, 12pt]{article}
   \usepackage{multirow}
          \usepackage{mathtools}

\usepackage[utf8]{inputenc}

\usepackage{amsmath}
\usepackage{amsthm}
\usepackage{amsfonts}
\usepackage{amssymb}
\usepackage{mathpazo}
\usepackage{float}
\usepackage{enumitem}

\usepackage{xcolor}

\newcommand{\R}{\mathbb{R}}

\renewcommand{\epsilon}{\varepsilon}

\newcommand{\A}{\mathbf{A}} 

\newcommand{\intersect}{\cap}
\newcommand{\union}{\cup}
\newcommand{\Intersect}{\bigcap}
\newcommand{\Union}{\bigcup}

\newcommand{\indicatorFunction}{\mathbb{1}}

\newcommand{\Product}{\prod}
\newcommand{\Sum}{\sum}

\newcommand{\freeParameter}{s}
\newcommand{\initialSubset}[1]{[#1]}

\numberwithin{equation}{section}
\newtheorem{theorem}{Theorem}
\newtheorem{lemma}[theorem]{Lemma}
\newtheorem{proposition}[theorem]{Proposition}
\newtheorem{corollary}[theorem]{Corollary}

\theoremstyle{definition}

\newtheorem{assumption}[theorem]{Assumption}
\theoremstyle{remark}

\newtheorem{remark}[theorem]{Remark}
\newtheorem{example}[theorem]{Example}

\usepackage{lscape}

\newcommand{\invariantForSup}{m}
\newcommand{\invariantForInf}{p}
\renewcommand{\a}{\mathbf{a}}

\begin{document}







\providecommand{\keywords}[1]
{
  \small	
  \textit{Keywords:   } #1
}

\providecommand{\msccode}[1]
{
  \small	
  \textit{2020 MSC:   } #1
}

\title{Probability bounds for $n$ random events under $(n - 1)$-wise independence}
\author{Karthik Natarajan $^{a,}$\footnote{Corresponding author: \texttt{karthik\char`_natarajan@sutd.edu.sg}} 
				\and Arjun Kodagehalli Ramachandra$^{b}$
				\and  Colin Tan$^{c}$ \\
        \small $^{a}$Engineering Systems and Design, Singapore University of Technology and Design \\
        \small $^{b}$Engineering Systems and Design, Singapore University of Technology and Design \\
	\small $^{c}$Engineering Systems and Design, Singapore University of Technology and Design \\
}

\date{03 November 2022}






\maketitle

\begin{abstract}
A collection of $n$ random events is said to be $(n - 1)$-wise independent if any $n - 1$ events among them are mutually independent. We characterise all probability measures with respect to which $n$ random events are $(n - 1)$-wise independent. We provide sharp upper and lower bounds on the probability that at least $k$ out of $n$ events with given marginal probabilities occur over these probability measures. The bounds are shown to be computable in polynomial time.
\end{abstract}

\noindent \keywords{$(n - 1)$-wise independence, mutual independence, probability bounds, polynomial time, Bonferroni bounds, probabilistic method}
%

\noindent \msccode{60E05}



\section{Introduction} \label{sec:intro}

Let $\Omega = \{\omega_J : \, J \subseteq [n]\}$
be the sample space freely generated by $n$ random events $A_1, \dots, A_n$, so that $\{\omega_J\} = \Intersect_{j \in J} A_j \intersect \Intersect_{j \notin J} \overline{A}_j =: \A^J$
			for all subsets $J \subseteq [n]$.
Here, as usual, for any integer $n \ge 0$, let $[n] := \{1, \dots, n\}$, let $\subseteq$ (resp.\ $\subset$) denote the subset (resp.\ proper subset) relation, let $\overline{A}$ denote the complement of $A$ and let $|J|$ denote the cardinality of a set $J$.
With $\Sigma$ as the $\sigma$-algebra of all subsets of $\Omega$,
the unique probability measure $P$ on the measurable space $(\Omega, \Sigma)$
    with respect to which $A_1, \dots, A_n$ are mutually independent
        and whose marginal probabilities are $P(A_j) =: a_j$ for all $j \in [n]$
	is given by:
\begin{equation} \label{eq: mutualIndependenceSolution}
P(\A^J) = \Product_{j \in J} a_j \times \Product_{j \notin J} (1 - a_j) =: \a^J,
	\quad \text{for all } J \subseteq [n].
\end{equation}				
We may relax the condition that the $n$ events are mutually independent
    to require only that every $(n - 1)$ events among $A_1, \dots, A_n$
        are mutually independent.
This weaker condition is sometimes known as \emph{$(n - 1)$-wise independence}. It is known that $(n - 1)$-wise independence does not imply that the $n$ events are mutually independent in general (see \cite{Wang,wangstoy} for counterexamples), although the converse is true. Comparisons of various notions of independence and dependence for random variables can be found in the literature \cite{Feller, Mukhopadhyay, Stoyanov, Wong}. The special case of Bernoulli random variables correspond to the setting of random events in this paper.

Bernstein \cite[p.\ 126]{Feller} constructed his classic example of $n = 3$ pairwise independent events which are not mutually independent. The example in \cite{Wang} generalizes his construction to $n \ge 3$ random events that are $(n - 1)$-wise independent but not mutually independent. We discuss this construction next.

\begin{example} [Construction of $(n-1)$-wise independent events] \label{ex:wang}
Let $a_j = 1/2$ for all $j \in [n]$ where $A_1,\ldots,A_{n-1}$ are mutually independent events and suppose that the event $A_n$ occurs given that an even number of events among $A_1,\ldots,A_{n-1}$ occur. It can be verified that these $n$ events are $(n - 1)$-wise independent but not mutually independent.

The induced probability measure is not the unique one with respect to which $A_1,\ldots,A_{n}$ are $(n - 1)$-wise independent but not mutually independent. For example, if we suppose instead that $A_n$ occurs given that an odd number of events among $A_1,\ldots,A_{n-1}$ occur, we obtain a distinct probability measure.
\end{example}
\subsection*{Overview}
In Section \ref{sec:characterization_(n-1)independence}, we characterize all probability measures with respect to which $n$ random events are $(n - 1)$-wise independent (Theorem \ref{thm: characterisationOfProbabilityMeasures}). Now we briefly outline the steps leading to Theorem \ref{thm: characterisationOfProbabilityMeasures}. First, in Proposition \ref{prop: formalEquivalence}, we identify a system of $2^n - 1$ equations linear in $P(\A^I)$ for all $I \subseteq [n]$ that is satisfied if and only if $A_1, \dots, A_n$ are $(n - 1)$-wise independent. We then relax the condition that $P(\A^I)$ are nonnegative and show in Lemma \ref{lem:1} and Corollary \ref{cor: characterisationOfUnsignedMeasures} that the solutions are parameterized by a single real parameter. Reimposing the condition that $P(\A^I)$ are nonnegative forces this parameter to lie within a compact interval which is identified in \eqref{eq: solveSystemForPositivityOfMeasure}.
Lemmas \ref{lem2}--\ref{lem: minimizeAtomicProbabilityForUpperBound} then explicitly identifies the endpoints of this interval in terms of the marginal probabilities. Theorem \ref{thm: characterisationOfProbabilityMeasures} then follows.

Section \ref{sec:atleastksharpbounds} applies Theorem \ref{thm: characterisationOfProbabilityMeasures} to identify sharp bounds on the probability that at least $k$ out of $n$ events which are $(n - 1)$-wise independent occurs (Theorem \ref{thm:atleastksharpbounds}) while Theorem \ref{thm:complexity} shows that these bounds are computable in polynomial time.

Examples \ref{ex:reductiontopairwise} and \ref{ex:reductiontobonferroni} in Section \ref{sec:examples} give cases when the newly derived bounds are instances of known universal bounds such as the classical Bonferroni bounds. Example \ref{ex:lovaszlocallemma} illustrates the connection of the results to the probabilistic method which provides conditions for the non-occurrence of ``bad" events when events are mostly independent. Example \ref{ex:comparisonmakarov} illustrates the usefulness of the bounds in providing robust estimates when mutual independence breaks down or when existing bounds are not sharp.


\section{Characterization of $(n-1)$-wise independence}\label{sec:characterization_(n-1)independence}
\begin{proposition} \label{prop: formalEquivalence}
A collection of $n$ random events $A_1, \dots, A_n$ is $(n - 1)$-wise independent with $P(A_j) = a_j$ for all $j \in [n]$
if and only if:
\begin{equation} \label{eq: n-1wiseIndependence}
\Sum_{I \supseteq J} P(\A^I) = \Product_{j \in J} a_j,
	\quad \text{for  all } J \subset [n].
\end{equation}
\end{proposition}

\begin{proof}
Since the $\A^I$'s are mutually exclusive, we have $\Sum_{I \supseteq J} P(\A^I) = P(\Intersect_{j \in J} A_j)
$ for any $J$.
Thus \eqref{eq: n-1wiseIndependence} is equivalent to
$$P(\Intersect_{j \in J} A_j) = \Product_{j \in J} a_j,
	\quad \text{for  all } J \subset [n],$$
which is in turn equivalent to the $(n - 1)$-wise independence of $A_1, \dots, A_n$.
\end{proof}

The system of linear equations in \eqref{eq: n-1wiseIndependence} is inhomogeneous with \eqref{eq: mutualIndependenceSolution}
	as a particular solution.
	
\begin{lemma} \label{lem:1}
The general solution of the system of linear equations:
\begin{equation} \label{eq: associatedHomogeousSystem}
\Sum_{I \supset J} P(\A^I) = 0,
	\quad \text{for  all } J \subset [n],
\end{equation}
in the $2^n$ variables $\{P(\A^J)\}_{J \subseteq [n]}$ (not necessarily nonnegative),
	is given by $P(\A^J) = (-1)^{|J|} s$, where $s \in \R$ is a free parameter.
\end{lemma}

\begin{proof}
The following inhomogeneous linear system has a unique solution, namely the mutually independent probability measure in \eqref{eq: mutualIndependenceSolution}:
$$\Sum_{I \supseteq J} P(\A^I) = \Product_{j \in J} a_j,
	\quad \text{for  all } J \subseteq [n].$$
 Hence the square coefficient matrix of its associated homogeneous linear system is invertible.
Thus, the equations of the homogeneous subsystem \eqref{eq: associatedHomogeousSystem} are linearly independent.
Therefore its solution space has dimension $1$,
    because there are $2^n - 1$ equations and $2^n$ variables. But, substituting $P(\A^I) = (-1)^{|I|} s$ into \eqref{eq: associatedHomogeousSystem}, where $s \in \R$ is a parameter, we have for any $J \subset [n]$:
\begin{align*}
\Sum_{I \supset J} (-1)^{|I|} s &= s \Sum_{q = |J|}^n \binom{n - |J|}{q - |J|}(-1)^q\\
&= s (-1)^{|J|} \Sum_{q = |J|}^n \binom{n - |J|}{q - |J|}(-1)^{q - |J|} \notag \\
&= s(-1)^{|J|}(1 - 1)^{n - |J|} \\
&= 0,
\end{align*}
where we use the Binomial Theorem in the third equality and the fact that $n - |J| > 0$ for all $J \subset [n]$ in the last equality.
Thus $P(\A^J) = (-1)^{|J|} s$.
\end{proof}

\noindent Recall that a measure $P$ is \emph{unitary} if $P(\Omega) = 1$.

\begin{corollary} \label{cor: characterisationOfUnsignedMeasures}
Every unitary measure $P$ (not necessarily nonnegative) on $(\Omega, \Sigma)$
	with $P(A_j) = a_j$ for all $j \in [n]$ and with respect to which $A_1, \dots, A_n$
    are $(n - 1)$-wise independent has the form:
\begin{equation} \label{eq:characterisationOfN-1wiseUnsignedMeasures}
P(\A^J) = \a^J + (-1)^{|J|} s, \quad \text{for all } J \subseteq [n],
\end{equation}		
for some scalar parameter $s \in \R$.
\end{corollary}

\begin{proof}
Add the particular solution $P(\A^J) = \a^J$
	(when the collection of events $\{A_1, \dots, A_n\}$ is mutually independent)
with the general solution of the associated homogenous linear system
	given in the previous lemma.
\end{proof}

\noindent To characterise when $P$  in \eqref{eq:characterisationOfN-1wiseUnsignedMeasures} is a valid probability measure, we have to ensure that the following $2^n$ nonnegativity conditions are satisfied:
\begin{equation} \label{eq: systemOfLinearInequalitiesToEnsurePositivity}
P(\A^J)=\a^J + (-1)^{|J|} s \ge 0, \quad \text{for  all } J \subseteq [n].
\end{equation}
Simplifying, this system gives:
\begin{equation}
\begin{cases}
 s \ge -\a^J, &\text{ for  all }J \subseteq [n]:\, |J| \text{ is even}, \\
s \le \a^J, &\text{ for  all } J \subseteq [n]:\, |J| \text{ is odd}.
\end{cases}
\end{equation}
Therefore it is a valid probability measure for all values of $s$ that satisfy:
\begin{equation} \label{eq: solveSystemForPositivityOfMeasure}
- \min_{J \subseteq [n]:\, |J| \text{ is even}} \a^J \le s \le  \min_{J \subseteq [n]:\, |J| \text{ is odd}} \a^J.
\end{equation}

\noindent From this point onwards, we make the following assumption.
\begin{assumption} \label{ass:order}
The events are ordered by nondecreasing value of their marginal probabilities, i.e. $a_1 \le \cdots \le a_n$.
\end{assumption}
\noindent The next lemma provides a lower bound on $\a^J$ for any set $J \subseteq [n]$, which will be used to establish the precise interval for the parameter $s$ in \eqref{eq: solveSystemForPositivityOfMeasure} and thus to identify all probability measures with respect to which $A_1, \dots, A_n$ are $(n - 1)$-wise independent.

\begin{lemma} \label{lem2}
$\a^J \ge \a^{\initialSubset{|J|}}$ for all $J \subseteq [n]$.
\end{lemma}

\begin{proof}
Use the notation $J \succeq I$ to denote $\a^J \ge \a^I$.
Say $|J| = \ell$. We need to show that $J \succeq \initialSubset{\ell} = \{1,\ldots,\ell\}$. If  $J=[\ell]$, we are done.
Otherwise, if $J \neq \initialSubset{\ell}$,
then there is a smallest index $r \le \ell$ that is not in $J$. 
Indeed, if there did not exist such an $r$, then the smallest index not in $J$ is strictly greater than $\ell$, hence $J \supset [\ell]$, which contradicts $|J| = \ell$. Hence the smallest index not in $J$, which we denote by $k$, satisfies $k > \ell$.

Let $J^\prime := J \union \{r\} \setminus \{k\}$
	be the set obtained by replacing $k$ with $r$ in $J$;	
		hence $|J^\prime| = |J| = \ell$ have the same cardinality.
Now, a common factor of $\a^J$ and $\a^{J^\prime}$ is
$C := \Product_{j \in J \setminus \{k\}} a_j
	\times \Product_{j \notin J \cup \{r\}} (1 - a_j)$,
	so
\begin{equation*}
\a^{J} - \a^{J^\prime} = C \big((1 - a_{r}) a_k - a_{r}(1- a_k) \big)
= C (a_k - a_{r}) \ge 0,
\end{equation*}		
since $a_k \ge a_{r}$ because $k > r$. Thus
\begin{equation*}
J \succeq J^\prime.
\end{equation*}
Now, if  $J^\prime=[\ell]$, we are done, else, repeating this procedure, we get a finite sequence of subsets, each of cardinality $\ell$ which terminates at $J^{\prime \dots \prime}=[\ell]$ after $|J \setminus [\ell]|$ replacements (since each iteration replaces exactly one element of $J \setminus [\ell]$ with one element of $[\ell] \setminus J$).
\end{proof}
\noindent We next define two integer invariants $p$ and $m$ of the ordered sequence of marginal probabilities.
These invariants are used to formulate the sharp lower bound on $\a^J$ for odd $|J|$ and that for even $|J|$ respectively.
Lemma \ref{lem2} is used to obtain these bounds.
Associate to $a_1 \le \cdots \le a_n$ an integer $\invariantForInf \in \{0, 1, \dots, \lfloor (n - 1)/2 \rfloor\}$
defined as:
\begin{multline} \label{eq: defineInvariantForInf}
p \mbox{ is the largest integer such that } \\
a_2 + a_3 \le \cdots \le a_{2 \invariantForInf} + a_{2 \invariantForInf + 1} \le 1.
\end{multline}

\begin{lemma} \label{lem: minimizeAtomicProbabilityForLowerBound}
If $|J|$ is odd, then:
\begin{equation}
\a^J \ge \a^{\initialSubset{2\invariantForInf + 1}},
\end{equation}
where $p$ is defined in \eqref{eq: defineInvariantForInf}.
\end{lemma}

\begin{proof}
Say $|J| = 2q + 1$.
First use Lemma \ref{lem2} to get $\a^J \ge \a^{\initialSubset{2q + 1}}$.
If $q = p$, we are done. Otherwise, either $q < p$ or $q > p$.
\paragraph{Case 1: Suppose $q < p$}
Then, we compare $\a^{\initialSubset{2q + 1}}$ with $\a^{\initialSubset{2q + 3}}$,
	which have $C := a_1 \cdots a_{2q + 1} (1 - a_{2q + 4}) \cdots (1 - a_n)$ as a common factor, hence
\begin{align*}
\a^{\initialSubset{2q + 1}} - \a^{\initialSubset{2q + 3}}
&= C\big(
(1 - a_{2q + 3})(1 - a_{2q + 2}) - a_{2q + 3} a_{2q + 2}
\big) \\
&= C\big(
1 - a_{2q + 2} - a_{2q + 3}
\big)\\
&\ge 0,
\end{align*}
since $a_{2q + 2} + a_{2q + 3} \le a_{2p } + a_{2p + 1} \le 1$.
Thus $\a^{\initialSubset{2q + 1}} \ge \a^{\initialSubset{2q + 3}}$. Repeating this process, we get
\begin{equation*}
\a^{\initialSubset{2q + 1}} \ge \a^{\initialSubset{2q + 3}} \ge \cdots \ge
\a^{\initialSubset{2p + 1}}.
\end{equation*}

\paragraph{Case 2: Suppose $q > p$} One can similarly show that $\a^{\initialSubset{2q + 1}} \ge a^{\initialSubset{2q - 1}}$ using $a_{2q} + a_{2q + 1} > 1$ because $p$ is the largest integer such that $a_{2p} + a_{2p + 1} \le 1$.
Repeat the process to get
$\a^{\initialSubset{2q + 1}} \ge \a^{\initialSubset{2q - 1}}
\ge \cdots \ge \a^{\initialSubset{2p + 1}}$.
\end{proof}
\noindent
Associate also to $a_1 \le \cdots \le a_n$ an integer $\invariantForSup \in \{0, 1, \dots, \lfloor n/2\rfloor\}$ defined as:
\begin{multline} \label{eq: defineInvariantForSup}
m \text{ is the largest integer such that } \\
a_1 + a_2 \le \cdots \le a_{2m - 1} + a_{2m} \le 1.
\end{multline}
The proof of the next lemma is similar to the previous lemma
and we omit it.

\begin{lemma} \label{lem: minimizeAtomicProbabilityForUpperBound}
If $|J|$ is even, then:
\begin{equation}
\a^J \ge \a^{\initialSubset{2\invariantForSup}},
\end{equation}
where $m$ is defined in \eqref{eq: defineInvariantForSup}.
\end{lemma}
\noindent This brings us to the following theorem.
\begin{theorem} \label{thm: characterisationOfProbabilityMeasures}
Let $p$ and $m$ be defined as in \eqref{eq: defineInvariantForInf} and \eqref{eq: defineInvariantForSup} respectively where Assumption \ref{ass:order} holds. Then every probability measure $P$ on $(\Omega, \Sigma)$
	with $P(A_j) = a_j$ for all $j \in [n]$ and with respect to which $A_1, \dots, A_n$ are $(n - 1)$-wise independent has the form:
$$
P(\A^J) = \a^J + (-1)^{|J|} s,\quad \mbox{for all} \;J \subseteq [n],$$ where $s$ is a scalar parameter satisfying:
\begin{equation} \label{eq: rangeOfFreeParameterToGetPositivity}
-\prod_{i=1}^{2m}a_i\prod_{i=2m+1}^{n}(1-a_i) \leq s \leq
\prod_{i=1}^{2p+1}a_i\prod_{i=2p+2}^{n}(1-a_i).
\end{equation}
\end{theorem}

\begin{proof}
By Corollary \ref{cor: characterisationOfUnsignedMeasures},
the conditions that
$P(A_j) = a_j$ for all $j \in [n]$
and $A_1, \dots, A_n$ are $(n - 1)$-wise independent with respect to $P$ entails that		
$
P(\A^J) = \a^J + (-1)^{|J|} s
$ for $J \subseteq [n]$ for some $s \in \R$. In order for $P$ to be a valid probability measure,
$s$ has to satisfy \eqref{eq: solveSystemForPositivityOfMeasure}. This gives:
\begin{align}
s &\in [-\min_{J \subseteq [n]: |J| \text{ is even}} \a^J ,  \min_{J \subseteq [n]: |J| \text{ is odd}} \a^J] \notag\\
&= [-\a^{\initialSubset{2m}}, \a^{\initialSubset{2p + 1}}],\notag
\end{align}
where the equality follows from Lemma \ref{lem: minimizeAtomicProbabilityForLowerBound} and Lemma \ref{lem: minimizeAtomicProbabilityForUpperBound}.
\end{proof}

\noindent
Previous results in the literature are limited to the construction of specific counterexamples showing that $(n - 1)$-wise independence does not imply the mutual independence of $n$ events (see \cite{Wang,wangstoy}). Theorem \ref{thm: characterisationOfProbabilityMeasures} comprehensively characterizes all such counterexamples.
Note that $s = 0$ corresponds to mutual independence.

\begin{remark}
Let us revisit the constructions given in Example \ref{ex:wang} in view of Theorem \ref{thm: characterisationOfProbabilityMeasures}. The probability measure where $A_n$ occurs given that an even number of events in $A_1, \dots, A_{n - 1}$ occur is given by
$$
P(\A^J) = \frac{1}{2^n} -  \frac{(-1)^{|J|}}{2^n},  \quad \text{for all } J \subseteq [n].$$
Since $s = - 1/2^n \in [-1/2^n, 1/2^n]$, it follows from Theorem \ref{thm: characterisationOfProbabilityMeasures}
that the events are $(n - 1)$-wise independent. Note that the events are not mutually independent since $s \neq 0$.

The other construction where $A_n$ occurs given that an odd number of events in $A_1, \dots, A_{n - 1}$ occur is the case of $s = 1/2^n$.
\end{remark}

\begin{proposition} \label{prop1}
Either $m = p$ or $m = p + 1$.
\end{proposition}

\begin{proof}
Let $k$ be the largest integer such that $a_i + a_{i + 1} \le 1$ for all $i \in [k]$.
If $k$ is odd, then $k = 2m - 1$ and $k - 1 = 2 p$,
hence $m = (k + 1)/2 = p + 1$.
On the other hand, if $k$ is even, then $k = 2 p $ and $k - 1 = 2 m - 1$, hence $m = k/2 = p$.
\end{proof}

\section{Probability bounds on at least $k$ events occurring}	\label{sec:atleastksharpbounds}		
\noindent In this section, we derive sharp bounds on the probability that at least $k$ out of $n$ events that are $(n - 1)$-wise independent occur. Bounds of this type under differing assumptions on the dependence structure of the random events have been studied (see \cite{rusch,boros1989}). Here we provide results for $(n-1)$-wise independence. From Theorem \ref{thm: characterisationOfProbabilityMeasures}, we can represent all such probabilities by:
\begin{align*}
P_{s}(n,k,\a)
&:= \sum_{q \ge k} \sum_{J \subseteq [n]: \, |J|=q} P(\A^J)\\
& = \sum_{q \ge k}^n \sum_{J \subseteq [n]: \, |J|=q} \left(\a^J + (-1)^{|J|} s\right),
\end{align*}
where $s \in [- \a^{\initialSubset{2m}}, \a^{\initialSubset{2p + 1}}]$ .
Then $P_0(n,k,\a)$ is the probability of occurrence of at least $k$ out of $n$ mutually independent events
with the given marginal probabilities $\a$.
We show that $P_{s}(n,k,\a)$ is linear in $s$.
Recall the binomial coefficient given by $\binom{z}{m} = z(z - 1) \cdots (z - m + 1)/m!$ for integers $z \ge 0$ and $n$.
\begin{lemma}\label{lem:atleastkevents}
For any integer $k \ge 0$,
\begin{align}\label{eq:atleastkintermsofalln}
P_{s}(n,k,\a) = P_0(n,k,\a) +(-1)^k\dbinom{n-1}{k-1} s,
\end{align}
where $s \in [- \a^{\initialSubset{2m}}, \a^{\initialSubset{2p + 1}}]$.
\end{lemma}
\begin{proof}
We have:
\begin{align*}
P_{s}(n,k,\a)
&= P_0(n, k, \a) + s \sum_{q \ge k} (-1)^q \sum_{J \subseteq [n] : \, |J| = q} 1 \\
&= P_0(n, k, \a) + s \sum_{q \ge k} (-1)^q \binom{n}{q}
\end{align*}	
The result then follows from the combinatorial identity $\sum_{q \ge k} (-1)^{q} \binom{n}{q} =  (-1)^k\binom{n-1}{k-1} $.
\end{proof}
We next derive sharp upper and lower bounds on the probability that at least $k$ out of $n$ events occur under $(n-1)$-wise independence.
\begin{theorem}\label{thm:atleastksharpbounds}
Let $p$ and $m$ be defined as in \eqref{eq: defineInvariantForInf} and \eqref{eq: defineInvariantForSup} respectively where Assumption \ref{ass:order} holds. Then every probability measure $P$ on $(\Omega, \Sigma)$
	with $P(A_j) = a_j$ for all $j \in [n]$ and with respect to which $A_1, \dots, A_n$ are $(n - 1)$-wise independent satisfies:
\begin{enumerate}[label=\roman*),wide=0pt]
\item For $k$ odd:
\begin{align}
P(n,k,\a)  &\ge P_0(n,k,\a) - \binom{n-1}{k-1}\prod_{i=1}^{2p+1} a_i \prod_{i=2p+2}^{n} (1-a_i),
\label{eq: lowerBoundForatleastkProbability2} \\
P(n,k,\a) &\le P_0(n,k,\a) + \binom{n-1}{k-1}\prod_{i=1}^{2m} a_i \prod_{i=2m+1}^{n} (1-a_i), \label{eq: upperBoundForatleastkProbability2}
\end{align}
\item For $k$ even:
\begin{align}
P(n,k,\a)  &\ge P_0(n,k,\a) - \dbinom{n-1}{k-1}\prod_{i=1}^{2m} a_i \prod_{i=2m+1}^{n} (1-a_i),
\label{eq: lowerBoundForatleastkProbability1} \\
P(n,k,\a)  &\le P_0(n,k,\a) + \dbinom{n-1}{k-1}\prod_{i=1}^{2p+1} a_i \prod_{i=2p+2}^{n} (1-a_i). \label{eq: upperBoundForatleastkProbability1}
\end{align}
\end{enumerate}
Moreover, all the bounds are sharp.
The lower bound for odd $k$ in \eqref{eq: lowerBoundForatleastkProbability2} and the upper bound for even $k$ in \eqref{eq: upperBoundForatleastkProbability1} is uniquely achieved with $P = P_{\a^{\initialSubset{2p + 1}}}$.
The upper bound for odd $k$ is \eqref{eq: upperBoundForatleastkProbability2} and the lower bound for even $k$ in \eqref{eq: lowerBoundForatleastkProbability1} is uniquely achieved with $P = P_{-\a^{\initialSubset{2m}}}$.
\end{theorem}
\begin{proof}
The result is obtained from Lemma \ref{lem:atleastkevents} and optimally selecting $\freeParameter$ in $[-\a^{\initialSubset{2m}},\a^{\initialSubset{2p + 1}}]$ from Theorem \ref{thm: characterisationOfProbabilityMeasures}.
\end{proof}

\noindent Probability bounds on the occurrence of at least $k$ out of $n$ events that are $\ell$-wise independent (i.e. every $\ell$ out of the $n$ events are mutually independent) have been studied for particular values of $\ell$. The case of $\ell = 1$ represents arbitrary dependence among the random events, for which the sharp upper bound is derived for $k = 1$ in \cite{boole} and for general $k$ in \cite{Ruger}.
At the other extreme is mutual independence ($\ell = n$), where the said probability is unique.
For $\ell = 2$, the sharp upper bound on the probability of the union ($k = 1$) of pairwise independent random events has been recently derived in \cite{ramanatarajan2021pairwise} and new bounds that are not necessarily sharp have been proposed for $k \ge 2$. Further, to the best of our knowledge, sharp bounds for other values of $\ell \in [3,n-1]$ have not been identified in the literature. Our results contribute to this line of work by finding sharp bounds for $l =n-1$.
We next demonstrate, as an immediate implication, the computability of the sharp bounds in Theorem \ref{thm:atleastksharpbounds}.

\begin{theorem} \label{thm:complexity}
The sharp upper and lower bounds in Theorem \ref{thm:atleastksharpbounds} are computable in polynomial time.
\end{theorem}
\begin{proof}
The value of $P_0(n,k,\a) $ is computable in polynomial time using dynamic programming. To see this, let $P_0(r,t,\a)$ denote the probability that at least $t$ events occur out of the first $r$ events where $r \geq t \geq 0$. Then the probabilities satisfy the recursive formula:
\begin{equation*}
 P_0(r,t,\a) = P_0(r-1,t-1,\a)p_r + P_0(r-1,t,\a)(1-p_r),
\end{equation*}
where the boundary conditions are $P_0(r,0,\a) = 1$ for $r \geq 0$ and $P_0(r,t,\a) = 0$ for $t > r$ (see \cite{hong2013poissonbinomial}).
The probability $P_0(n,k,\a) $ is hence computable in $O(n^2)$ time.  Since the additional term in the formulas \eqref{eq: lowerBoundForatleastkProbability2}-\eqref{eq: upperBoundForatleastkProbability1} is efficiently computable using sorting, evaluating binomial coefficients and multiplication, all the bounds are computable in polynomial time; specifically $O(n^2)$ time.
\end{proof}

\noindent The sharp bounds for $k = 1$ (union of events) and $k = n$ (intersection of events) are detailed next.
\begin{corollary} \label{unionintersect}
Let $p$ and $m$ be defined as in \eqref{eq: defineInvariantForInf} and \eqref{eq: defineInvariantForSup} respectively where Assumption \ref{ass:order} holds. Then every probability measure $P$ on $(\Omega, \Sigma)$
	with $P(A_j) = a_j$ for all $j \in [n]$ and with respect to which $A_1, \dots, A_n$ are $(n - 1)$-wise independent satisfies:
\begin{enumerate}[label=\roman*),wide=0pt]
\item For the union of events:
\begin{align}
P(\Union_{j = 1}^n A_j)  &\ge 1 - \left(\prod_{i=1}^{2p+1}(1-a_i)+\prod_{i=1}^{2p+1}a_i\right)\prod_{i=2p+2}^{n}(1 - a_i),
\label{eq: lowerBoundForUnionProbability} \\
P(\Union_{j = 1}^n A_j)  &\le  1 - \left(\prod_{i=1}^{2m}(1-a_i)-\prod_{i=1}^{2m}a_i\right)\prod_{i=2m+1}^{n}(1 - a_i), \label{eq: upperBoundForUnionProbability}
\end{align}
\item For the intersection of an even number of events:
\begin{align}
P(\Intersect_{j = 1}^n A_j) &\ge \prod_{i=1}^{2m}a_i\left(\prod_{i=2m+1}^{n}a_i-\prod_{i=2m+1}^{n}(1-a_i)\right),
\label{eq: lowerBoundForIntersectProbabilityeven} \\
P(\Intersect_{j = 1}^n A_j) &\le\prod_{i=1}^{2p+1}a_i\left(\prod_{i=2p+2}^{n}a_i+\prod_{i=2p+2}^{n}(1-a_i)\right),
\label{eq: upperBoundForIntersectProbabilityeven}
\end{align}
\item For the intersection of an odd number of events:
\begin{align}
P(\Intersect_{j = 1}^n A_j) &\ge  \prod_{i=1}^{2p+1}a_i\left(\prod_{i=2p+2}^{n}a_i-\prod_{i=2p+2}^{n}(1-a_i)\right),
\label{eq: lowerBoundForIntersectProbabilityodd} \\
P(\Intersect_{j = 1}^n A_j) &\le \prod_{i=1}^{2m}a_i\left(\prod_{i=2m+1}^{n}a_i+\prod_{i=2m+1}^{n}(1-a_i)\right)\label{eq: upperBoundForIntersectProbabilityodd}.
\end{align}
\end{enumerate}
Each of these bounds is sharp and is achieved by a unique probability measure $P(\A^J) = \a^J + (-1)^{|J|}s$, where either $\freeParameter =-\a^{\initialSubset{2m}}$ or
$\freeParameter =\a^{\initialSubset{2p + 1}}$.
\end{corollary}

\begin{proof}
With $k=1$, we have $P_0(n,1,\a)=1-\prod_{i=1}^n (1-a_i)$ and the result immediately follows from \eqref{eq: lowerBoundForatleastkProbability2}  and \eqref{eq: upperBoundForatleastkProbability2} in Theorem \ref{thm:atleastksharpbounds}. With $k=n$, we have $P_0(n,n,\a)=\prod_{i=1}^n a_i$ and the result immediately follows from Theorem \ref{thm:atleastksharpbounds} depending on whether $k = n$ is even or odd.
\end{proof}

\section{Examples}\label{sec:examples}
\noindent In this section, we discuss several examples to illustrate the connection of the newly proposed bounds with existing bounds and provide numerical evidence of the quality of the bounds.
\begin{example} [Bounds for $n = 3$ pairwise independent events]\label{ex:reductiontopairwise}
For $n=3$ events, $(n-1)$-wise independence is pairwise independence.
In this case from \eqref{eq: defineInvariantForSup} and \eqref{eq: upperBoundForUnionProbability} where $a_1 \leq a_2 \leq a_3$, we obtain the sharp upper bound on the union as:
\begin{equation} \label{eq: sharpUpperBoundForThreeEvents}
P(\Union_{j = 1}^3 A_j)
 \le  \min\left(a_1+a_2+a_3-a_3(a_1+a_2),1\right)
\end{equation}
Another proof of the sharpness was given in \cite{ramanatarajan2021pairwise}.
Kounias \cite{kounias} showed that every probability measure
satisfies $P(\Union_{j = 1}^3 A_j) \le \min\left(a_1+a_2+a_3-a_{31}+a_{32}),1\right)$, where $a_{ij} := P(A_i \intersect A_j)$ for $i, j \in [3]$ denotes the bivariate joint probability.
Therefore, \eqref{eq: sharpUpperBoundForThreeEvents} entails that the upper bound of Kounias is achieved by some probability measure with respect to which $A_1, A_2, A_3$ are pairwise independent.

For the sharp lower bound, from \eqref{eq: defineInvariantForInf} and \eqref{eq: lowerBoundForUnionProbability}, we get:
\begin{multline*}
P(\Union_{j = 1}^3 A_j)
\ge \max(a_2+a_3-a_2a_3, \\
a_1+a_2+a_3-a_1a_2-a_1a_3-a_2a_3),
\end{multline*}
Similarly, a corresponding universal lower bound of Kounias in terms of bivariate joint probabilities is therefore achievable under pairwise independence.

Likewise for the intersection of three pairwise independent events, we can verify that the sharp bounds are given as:
\begin{equation*}
P(\Intersect_{j = 1}^3 A_j) \le  \min(a_1a_2,
(1 - a_1)(1 - a_2)(1 - a_3) + a_1 a_2 a_3),
\end{equation*}
and
\begin{equation*}
P(\Intersect_{j = 1}^3 A_j) \ge \max\left(a_1(a_2+a_3-1),0\right).
\end{equation*}
An alternative proof of the sharpness of the lower bound is given in \cite{ramanatarajan2021pairwise}.
\end{example}

\begin{example} [Bonferroni bounds]\label{ex:reductiontobonferroni}
Suppose the sum of the two largest marginal probabilities satisfies $a_{n-1} + a_{n} \leq 1$ and $n$ is even. Then $m=n/2$ in \eqref{eq: defineInvariantForSup}, hence from \eqref{eq: upperBoundForUnionProbability}, we get the sharp upper bound on the union:
\begin{align*}
P(\Union_{j = 1}^n A_j) &\le 1 - \prod_{i=1}^{n}(1 - a_i) + \prod_{i=1}^{n}a_i \\
&= \sum_{k = 0}^{n - 2} (-1)^k \sum_{1 \le i_0 < \cdots < i_k \le n} a_{i_0} \cdots a_{i_k}.
\end{align*}
Bonferroni \cite{Bonferroni} showed that every probability measure
satisfies
\begin{equation}
P(\Union_{j = 1}^n A_j) \le \sum_{k = 0}^{n - 2} (-1)^k \sum_{1 \le i_0 < \cdots < i_k \le n} a_{i_0 \cdots i_k},
\end{equation}
where $a_{i_0 \cdots i_k} := P(A_{i_0} \intersect \cdots \intersect A_{i_k})$ for $i_0, \dots, i_k \in [n]$ is the joint probability.
Therefore the Bonferroni upper bound is achieved by some probability measure with respect to which $A_1, \dots, A_n$ are $(n - 1)$-wise independent, in this case.

Similarly if $a_{n - 1} + a_n \le 1$ and $n$ is odd,
then $p = (n - 1)/2$ in \eqref{eq: defineInvariantForInf} and thus the sharp lower bound in \eqref{eq: lowerBoundForUnionProbability} becomes:
\begin{align*}
P(\Union_{j = 1}^n A_j) &\ge 1 - \prod_{i=1}^{n}(1 - a_i) - \prod_{i=1}^{n}a_i \\
&= \sum_{k = 0}^{n - 2} (-1)^k \sum_{1 \le i_0 < \cdots < i_k \le n} a_{i_0} \cdots a_{i_k}.
\end{align*}
Again, a corresponding lower bound of Bonferroni in terms of joint probabilities of up to $n - 1$ events is thus achievable under $(n - 1)$-wise independence.
\end{example}

The next example shows the connection of the bound to the probabilistic method which has proved to be very useful tool in combinatorics (see \cite{alon}).
\begin{example} [Probabilistic method]\label{ex:lovaszlocallemma} Suppose there are $n$ random ``bad'' events, each of which occurs with probability $a_j$ for $j \in [n]$. When the events are mutually independent, the probability of no bad event occurring is strictly positive when the probability of each bad event is strictly less than 1 (namely $\max_j  a_j < 1$). 
On the other hand, if the events can be arbitrarily dependent, from Boole's union bound \cite{boole}, the sum of the probabilities must be strictly less than 1 (namely $\sum_j  a_j < 1$) to guarantee the same. The Lov\'{a}sz local lemma \cite{erdos} is a powerful tool that allows one to relax the assumption of mutual independence to weak dependence while allowing for the probability of each bad event to be fairly large and still guarantee that no bad event occurs with strictly positive probability. Specifically consider a graph $G$ on $n$ nodes where each node $i \in [n]$ is associated with an event $A_i$ and $A_i$ is independent of the collection of events $\{A_j: (i,j) \notin G\}$ for each $i \in [n]$. If $G$ has maximum degree $d$ and $\max_i a_i \leq 1/4d$, then the probability of no bad event occurring satisfies (see \cite{erdos,tetali}):
\begin{align}
P(\Intersect_{j = 1}^n \overline{A}_j)
\ge \prod_{i=1}^{n}(1-2a_i) > 0.\label{eq:lovaszlower}
\end{align}
Computing the tightest lower bound in terms of the dependency graph is known to be NP-complete \cite{shearer}. More generally, in \cite{shearer} it was shown that for $d \geq 2$, $\max_i a_i < (d-1)^{d-1}/d^d$ and for $d = 1$, $\max_i a_i < 1/2$ guarantees that there is a strictly positive probability that no bad event occurs. For the specific case of $d=1$, we can compare our results with the lower bound as shown next (although the Lov\'{a}sz local lemma holds more generally for lesser independence with $d \ge2 $). When the events are $(n-1)$-wise independent, using \eqref{eq: upperBoundForUnionProbability}, the probability that none of the events occur is strictly positive if $a_n< 1$ and $a_1+a_2 < 1$.
Indeed, then $(1 - a_1)(1 - a_2) > a_1 a_2$ and $(1 - a_{2k - 1})(1 - a_{2k}) \ge a_{2k - 1} a_{2k}$ for all $k \in \{2, \dots, m\}$. Hence $\prod_{i=1}^{2m}(1-a_i) = \prod_{k = 1}^{m}((1 - a_{2k - 1})(1 - a_{2k})) > \prod_{k = 1}^{m} (a_{2k - 1} a_{2k}) = \prod_{i=1}^{2m}a_i$ and from \eqref{eq: upperBoundForUnionProbability}:
$$
P(\Intersect_{j = 1}^n \overline{A}_j)
\ge \left(\prod_{i=1}^{2m}(1-a_i)-\prod_{i=1}^{2m}a_i\right)\prod_{i=2m+1}^{n}(1 - a_i)
> 0.
$$

\noindent When all the marginal probabilities $a_1=\ldots=a_n= a$, are identical, the condition $a_1 + a_2 < 1$ gives $a < 1/2$ which exactly corresponds to the condition identified in \cite{spencer,shearer} for $d = 1$. It is easy to verify that the lower bound on the probability of no bad event occurring in this case is given by $(1-a)^n-a^n$ for $n$ even and $(1-a)^n-a^{n-1}(1-a)$ for $n$ odd which is the sharp lower bound instance wise. In comparison, the lower bound identified above in \eqref{eq:lovaszlower} is $(1-2a)^n$. For example with $n = 6$ and $a = 0.1$, the sharp lower bound is $0.53144$ while the weaker lower bound is $0.262144$.
In fact for $a = 1/2$, the first construction in Example \ref{ex:wang} has a zero probability that no bad event occurs since $A_n$ must occur when none of the events in $\{A_1,A_2,\ldots A_{n-1}\}$ occur.
\end{example}
\noindent We next provide a numerical example to illustrate the performance of the bounds in Theorem \ref{thm:atleastksharpbounds} and compare it with an existing bound.
Specifically tail probability bounds on the sum of two random variables given their marginal distribution functions were derived by Makarov in \cite{makarov1982}. We adopt these closed-form bounds also known as ``standard" bounds in our context as follows.
 Given that $n$ random events $A_1,\ldots A_n$ with respective marginal probabilities $a_1 \le \cdots \le a_n$ are $(n - 1)$-wise independent, define two random variables as follows: $Y_1=\sum_{i=1}^{n-1} \indicatorFunction_{A_i},\;Y_2=\indicatorFunction_{A_n}$
where $ \indicatorFunction_{A}$ is the indicator function of event $A$ occurring. Here $ Y_1 \sim \operatorname{PoissonBinomial}(n-1,a_1,a_2,\ldots a_{n-1})$ is an integer random variable taking values in $[0,n-1]$ while $Y_2 \sim \operatorname{Bernoulli}(a_n)$.
Let $F_1$ and $F_2$ be the resepective distribution functions of $Y_1$ and $Y_2$. Then the Makarov upper bound for the probability that the sum of $Y_1$ and $Y_2$ is at least an integer $k \in [n]$ is given from \cite{makarov1982,rusch} as follows:
\begin{equation}\label{eq: makarovupper}
\begin{array}{rll}
 P(Y_1+Y_2\ge k)& \leq
 & \min(2-(F_{1}\vee F_{2})^{-}(k),1),
\end{array}
\end{equation}
where $(F_{1}\vee F_{2})^{-}(k)=\underset{u \in \mathbb{R}}{\max}(F_{1}(k-u)^{-} +F_{2}(u))$ is the left continuous version of the supremum convolution $F_{1}\vee F_{2}$. Since $Y_2$ is a Bernoulli random variable, it is sufficient to maximize over $u \in \{0,1\}$ and thus we have :
 \begin{align*}
& (F_{1}\vee F_{2})^{-}(k)
 & =&\max(F_{1}(k-1) +a_n,\;F_{1}(k-2) +1).
  \end{align*}
 The Makarov lower bound can be similarly derived as
\begin{equation}\label{eq: makarovlower}
\begin{array}{rll}
 P(Y_1+Y_2\ge k)
\ge\max(1-\min(F_{1}(k),F_{1}(k-1) +a_n),0).
\end{array}
\end{equation}
We next illustrate through a numerical example that the Makarov bound is not sharp in general under $(n-1)$-wise independence since we lose out on using additional independence information available in our context. For example, our bounds assume that any $n-2$ events from the first $n-1$ events $A_1,\ldots,A_{n-1}$ along with the last event $A_n$ are mutually independent while the Makarov bounds do not assume so.
\begin{example} [Numerical example] \label{ex:comparisonmakarov} Here we compute the exact probability for $n = 8$ with identical marginal probabilities $a_i = a \in \{0.1,0.2,0.3,0.4,0.5\}$ for different values of $k$ assuming mutual independence. In addition we compute the sharp lower and upper bounds with $7$-wise independence from Theorem \ref{thm:atleastksharpbounds} (here $p = 3$ and $m = 4$ for all considered values of $a$). We also provide the Makarov lower and upper bounds from (\ref{eq: makarovlower}) and (\ref{eq: makarovupper}) to highlight that if more information is known on the independence of the random variables, we can exploit it tightening the bounds.
\begin{table}[H]
\footnotesize
\caption{$k=1$ to $k=4$ - For each value of $a$, the first row provides the Makarov lower bound from \eqref{eq: makarovlower}, the second row provides the sharp lower bound with $7$-wise independence, the third row provides the exact value with $8$ mutually independent events, the fourth row provides the sharp upper bound with $7$-wise independence and the fifth row provides the Makarov upper bound from \eqref{eq: makarovupper}} \label{tab:heteroall3vstight1}
\begin{center}
\scriptsize{\begin{tabular}
{|l|c|c|l|l|l|l|l|l|l|}
  \hline
\mbox{a} & $k = 1$ & $k = 2$ & $k = 3$ & $k = 4$  \\ \hline
0.1 & 4.6953e-01 & 8.6895e-02 & 5.0243e-03 & 4.3165e-04\\
0.1 & 5.6953e-01 & 1.8690e-01 & 3.8090e-02 & 5.0240e-03\\
0.1 & 5.6953e-01 & 1.8690e-01 & 3.8092e-02 & 5.0244e-03\\
0.1 & 5.6953e-01 & 1.8690e-01 & 3.8092e-02 & 5.0275e-03\\
0.1 & 1.0000e+00 & 5.6953e-01 & 1.8690e-01 & 3.8092e-02\\
 \hline
0.2 & 6.3223e-01 & 2.9668e-01 & 5.6282e-02 & 1.0406e-02\\
0.2 & 8.3222e-01 & 4.9667e-01 & 2.0287e-01 & 5.6192e-02\\
0.2 & 8.3223e-01 & 4.9668e-01 & 2.0308e-01 & 5.6282e-02\\
0.2 & 8.3223e-01 & 4.9676e-01 & 2.0314e-01 & 5.6640e-02\\
0.2 & 1.0000e+00 & 8.3223e-01 & 4.9668e-01 & 2.0308e-01\\
\hline
0.3 & 7.4470e-01 & 4.4823e-01 & 1.9410e-01 & 5.7968e-02\\
0.3 & 9.4220e-01 & 7.4424e-01 & 4.4501e-01 & 1.9181e-01\\
0.3 & 9.4235e-01 & 7.4470e-01 & 4.4823e-01 & 1.9410e-01\\
0.3 & 9.4242e-01 & 7.4577e-01 & 4.4960e-01 & 1.9946e-01\\
0.3 & 1.0000e+00 & 9.4235e-01 & 7.4470e-01 & 4.4823e-01\\
\hline
0.4 & 8.9362e-01 & 6.8461e-01 & 4.0591e-01 & 1.7367e-01\\
0.4 & 9.8222e-01 & 8.8904e-01 & 6.6396e-01 & 3.8298e-01\\
0.4 & 9.8320e-01 & 8.9362e-01 & 6.8461e-01 & 4.0591e-01\\
0.4 & 9.8386e-01 & 9.0051e-01 & 6.9837e-01 & 4.4032e-01\\
0.4 & 1.0000e+00& 9.8320e-01 & 8.9362e-01 & 6.8461e-01\\
\hline
0.5 & 9.6484e-01 & 8.5547e-01 & 6.3672e-01 & 3.6328e-01\\
0.5 & 9.9219e-01 & 9.3750e-01 & 7.7344e-01 & 5.0000e-01\\
0.5 & 9.9609e-01 & 9.6484e-01 & 8.5547e-01 & 6.3672e-01\\
0.5 & 1.0000e+00 & 9.9219e-01 & 9.3750e-01 & 7.7344e-01\\
0.5 & 1.0000e+00 & 9.9610e-01 & 9.6484e-01 & 8.5547e-01\\
\hline
\end{tabular}}
\end{center}
\end{table}

\begin{table}[H]
\footnotesize
\caption{$k=5$ to $k=8$} \label{tab:heteroall3vstight2}
\begin{center}
\scriptsize{\begin{tabular}[t]
{|l|c|c|l|l|l|l|l|l|l|}
  \hline
\mbox{a} & $k = 5$ & $k = 6$ & $k = 7$ & $k = 8$ \\ \hline
0.1 & 2.3410e-05 & 7.3000e-07 & 9.9999e-09 & 0.0000e+00\\
0.1 & 4.2850e-04 & 2.3200e-05 & 1.0000e-07 & 0.0000e+00\\
0.1 & 4.3165e-04 & 2.3410e-05 & 7.3000e-07 & 1.0000e-08\\
0.1 & 4.3200e-04 & 2.5300e-05 & 8.0000e-07 & 1.0000e-07\\
0.1 & 5.0244e-03 & 4.3165e-04 & 2.3410e-05 & 7.3000e-07\\
 \hline
0.2 & 1.2314e-03 & 8.4480e-05 & 2.5600e-06 & 0.0000e+00\\
0.2 & 1.0048e-02 & 1.1776e-03 & 1.2800e-05 & 0.0000e+00\\
0.2 & 1.0406e-02 & 1.2314e-03 & 8.4480e-05 & 2.5600e-06\\
0.2 & 1.0496e-02 & 1.4464e-03 & 1.0240e-04 & 1.2800e-05\\
0.2 & 5.6282e-02 & 1.0406e-02 & 1.2314e-03 & 8.4480e-05\\
\hline
0.3 & 1.1292e-02 & 1.2903e-03 & 6.5610e-05 & 0.0000e+00\\
0.3 & 5.2610e-02 & 9.9144e-03 & 2.1870e-04 & 0.0000e+00\\
0.3 & 5.7968e-02 & 1.1292e-02 & 1.2903e-03 & 6.5610e-05\\
0.3 & 6.0264e-02 & 1.4507e-02 & 1.7496e-03 & 2.1870e-04\\
0.3 & 1.9410e-01 & 5.7968e-02 & 1.1292e-02 & 1.2903e-03\\
\hline
0.4 & 4.9807e-02 & 8.5200e-03 & 6.5536e-04 & 0.0000e+00\\
0.4 & 1.3926e-01 & 3.6045e-02 & 1.6384e-03 & 0.0000e+00\\
0.4 & 1.7367e-01 & 4.9807e-02 & 8.5197e-03 & 6.5536e-04\\
0.4 & 1.9661e-01 & 7.0451e-02 & 1.3107e-02 & 1.6384e-03\\
0.4 & 4.0591e-01 & 1.7367e-01 & 4.9807e-02 & 8.5197e-03\\
\hline
0.5 & 1.4453e-01 & 3.5156e-02 & 3.9063e-03 & 0.0000e+00\\
0.5 & 2.2656e-01 & 6.2500e-02 & 7.8125e-03 & 0.0000e+00\\
0.5 & 3.6328e-01 & 1.4453e-01 & 3.5156e-02 & 3.9063e-03\\
0.5 & 5.0000e-01 & 2.2656e-01 & 6.2500e-02 & 7.8125e-03\\
0.5 & 6.3672e-01 & 3.6328e-01 & 1.4453e-01 & 3.5156e-02\\
\hline
\end{tabular}}
\end{center}
\end{table}
\noindent As it can be observed, the sharp bounds with $(n-1)$-wise independence clearly improve upon the Makarov bounds, especially as $k$ increases (for the same $a$) and $a$ decreases (for the same $k$), where the bounds can be a couple or more magnitude of orders apart. In other words, the sharp bounds especially provide value in the regime where the right tail probabilities are more constrained i.e. large $k$ and small $a$. Such bounds are useful in providing robust estimates of the probabilities when the assumption of mutual independence breaks down.
\end{example}

\subsection*{Acknowledgements}
\noindent The research of the first and third authors was partly supported by MOE Academic Research Fund Tier 2 grant T2MOE1906, ``Enhancing Robustness
770 of Networks to Dependence via Optimization''. The authors would like to thank the Area Editor Henry Lam, the Associate Editor and the anonymous reviewer for valuable comments.






\end{document}